\documentclass[12pt]{amsart}
\usepackage{amsmath,amssymb,amsthm,amscd}
\numberwithin{equation}{section} 

\def\R{{\mathfrak R}}

\newtheorem{prop}{Proposition}[section]
\newtheorem{theo}[prop]{Theorem}
\newtheorem{lemm}[prop]{Lemma}
\newtheorem{coro}[prop]{Corollary}
\newtheorem{rema}[prop]{Remark}
\newtheorem{exam}[prop]{Example}
\newtheorem{defi}[prop]{Definition}

\def\begeq{\begin{equation}}
\def\endeq{\end{equation}}
\def\p{\partial}
\def\cF{\mathcal F}
\def\lf{\left}
\def\ri{\right}
\def\e{\epsilon}
\def\ol{\overline}
\def\R{\Bbb R}
\def\wt{\widetilde}
\def\tOmega{\wt\Omega}
\def\tM{\wt M}
\begin{document}
\title{Quasi-local mass and the existence of horizons
}

\author{Yuguang Shi$^1$}
\address{Key Laboratory of Pure and Applied mathematics, School of Mathematics Science, Peking University,
Beijing, 100871, P.R. China.} \email{ygshi@math.pku.edu.cn}

\author{Luen-Fai Tam$^2$}

\thanks{$^1$Research partially supported by NSF grant of China}

\thanks{$^2$Research partially supported by Earmarked Grant of Hong
Kong \#CUHK403005}

\address{Department of Mathematics, The Chinese University of Hong Kong,
Shatin, Hong Kong, China.} \email{lftam@math.cuhk.edu.hk}

\renewcommand{\subjclassname}{%
  \textup{2000} Mathematics Subject Classification}
\subjclass[2000]{Primary 83C57  ; Secondary 53C44  ,  }

\date{November,  2005}

\begin{abstract}In this paper, we obtain lower bounds for the Brown-York
quasilocal mass and the Bartnik quasilocal mass for compact three
manifolds with smooth boundaries. As a consequence, we derive
 sufficient conditions for the existence of horizons for a certain
class of compact manifolds with boundary and some asymptotically
flat complete manifolds.
 The method is based on analyzing Hawking mass and
inverse mean curvature flow.

\end{abstract}

\maketitle \markboth{Yuguang Shi and Luen-Fai Tam} {Quasi-local
mass and the existence of horizons}

\section{Introduction}

In this work, we will discuss the relations between different
kinds of quasi-local mass and use them to derive   sufficient
conditions for the existence of horizons in the time symmetric
case. In the time symmetric case, a horizon is defined to be a
compact minimal surface which is the boundary of some open sets
(see \S3 for more precise definition).

 In 1972, Thorne made the following conjecture, which later became
known as the hoop conjecture(see \cite{F}): {\it Black holes with
horizons form when and only when a mass $M$ gets compacted into a
region whose circumference in every direction is $\mathcal{C}\leq
4 \pi M$.} The conjecture is loosely formulated . Several concepts
such as mass, circumference etc.   are not clearly defined. Hence,
this conjecture allows many different precise interpretations.

 In
1983, Schoen and Yau in \cite{SY}, and later Yau in \cite{Y}
define a   kind of radius of a bounded region, and derive an upper
bound for this radius for a region in  spacetime  without apparent
horizons in terms of the lower bound of  mass density. By this and
the study of the obstruction to the existence of regular solution
to the Jang equation, they obtain some important results for the
existence of black holes in the spirit of the hoop conjecture.

In the time symmetric case, the method does not work because the
Jang equation always has a solution. For this case, we will use
the inverse mean curvature flow of Huisken and Ilmanen \cite{HI}.
We will define a quantity $m(\Omega)$ on a compact three manifold
$\Omega$ with boundary which involves Hawking mass of some subsets
of $\Omega$ and some other geometric quantities  of $\Omega$, see
definition in (\ref{newmass}). Then we compare this with the
Brown-York mass $m_{BY}(\p\Omega)$ of $\p\Omega$ of   a simply
connected compact manifold with nonnegative scalar curvature and
with connected smooth boundary which has  positive Gauss curvature
and positive mean curvature with respect to the outward normal. It
is well-known that $\p\Omega$ can be isometrically embedded in
$\R^3$.  The first result is as follows: If $\Omega$ contain no
horizons, then $m_{BY}(\p\Omega)\ge m(\Omega)$. From this one can
prove that if $m(\Omega)\ge 2R$ where $R$ is the radius of the
smallest circumscribed ball of
 $\p\Omega $ in $\R^3$, then $\Omega$ must contains a horizon. In
 particular, if $m(\Omega)\ge 2\,\text{diam}(\p\Omega)$, then
 $\Omega$ contains a horizon. Examples satisfying these conditions
 are given.

 We will prove that if $\Omega$
 contains a {\it round sphere} $\p E$ in the sense of \cite{CY}
 such that its Hawking mass is larger than $m_{BY}(\p\Omega)$,
 then $\Omega$ also contains a horizon. However, it is unclear if one can
 find a round sphere satisfying the condition.

It turns out that the quasi-local mass $m_B(\Omega)$ (see the
definition in \S3) introduced by Bartnik for a compact manifold
$(\Omega,g)$ with smooth boundary with an admissible extension is
also bounded below by $m(\Omega)$. It was observed by Walter Simon
(see \cite{B3}) that if
 $ (\Omega, g)$   has an admissible extension and
 suppose $(\Omega,g)$ is isometrically embedded
 in an asymptotically flat and complete manifold $M$ with
nonnegative scalar curvature so that the ADM   mass $m_{ADM}(M)$
of $M$ is less than $m_B (\Omega)$,   then $M$ must contains a
horizon. Our lower bound for $m_B(\Omega)$ implies that if
$m_{ADM}(M)<m(\Omega)$, then $M$ must contain  a horizon.

We should emphasis that the quantity $m(\Omega)$ defined in (\ref
{newmass}) is nontrivial, in the sense that $m(\Omega)\ge0$ and is
zero only if it is locally flat. One can prove that it is actually
a domain in $\R^3$ in some cases.

The basic outline of the paper is as follows. In Section 1,  we
first construct examples to motivate the definition of
$m(\Omega)$. Then we prove   the positivity of $m(\Omega)$ (see
Theorem 1.1).

In Section 2, we compare  $m(\Omega)$, $m_{BY}(\Omega)$   to give
sufficient conditions for the existence of horizons for compact
manifold with boundary as mentioned above. We will give another
lower bound for the Brown-York mass. In particular, we prove that
the Hawking mass of the boundary is dominated by the Brown-York
mass of the domain.

In Section 3, we prove that   $m(\Omega)$ is bounded above by the
Bartnik mass $m_B (\Omega)$ for  $\Omega$ which has an admissible
extension. We will  also discuss some properties of Bartnik
quasi-local mass including the conjecture of Bartnik that
$m_B(\Omega)$ is realized by a  static admissible extension, see
\cite{B3} for more details of the conjecture.

The main analytical tool in this paper is the inverse mean
curvature flow which has been studied by   Huisken and   Ilmanen
\cite{HI}. We obtain our results by studying the obstruction for
the monotonicity of Hawking mass under this flow.

\section  {Hawking mass of subsets of a domain}

In this section, we will introduce a quantity involving
 Hawking mass of some subsets of a compact three manifold
  $(\Omega,g)$ with smooth boundary which will be used to give a
  condition for the existence
  of stable minimal spheres on a compact manifold with boundary.
   All Riemannian manifolds in this work are assumed to be oriented and
   connected with dimension three.

  To motivate the definition, let us construct some examples.

  \begin{prop}\label{examples} There exist asymptotically flat metrics $g_1$, $g_2$ with
  nonnegative scalar curvature on $\R^3$ such that $g_1=g_2$
  outside some compact set, $(\R^3,g_1)$ contains a stable minimal
  sphere but $(\R^3,g_2)$ does not contain any compact minimal
  surfaces.
  \end{prop}

  Recall that   an asymptotically flat (AF) three manifold
$(M,g)$ with one end   is a complete manifold
 with nonnegative scalar curvature which is in $L^1(M)$ such that for some compact
 set $K$ of $M$, $M\setminus K$ is diffeomorphic to $\R^3\setminus
 B_R(0)$ for some $R>0$ and in the standard coordinates in $\R^3$,
 the metric $g$ satisfies:

 \begin{equation}\label{AF1}
   g_{ij}=\delta_{ij}+b_{ij}
\end{equation}
with
\begin{equation}\label{AF2}
   ||b_{ij}||+r||\p b_{ij}||+r^2||\p\p b_{ij}||=O(r^{-1})
\end{equation}
 where $r$ and $\p$ denote the Euclidean distance and standard derivative
 operator on $\R^3$.

   Recall also that the ADM mass of $M$ is
 defined as

\begin{equation}\label{ADMmass}
m_{ADM}(M)=\lim_{r\to\infty}\frac1{16\pi}\int_{S_r}(  g_{ii,j}
-g_{ij,i})\nu^j
    d\sigma_r
\end{equation}
where $S_r$ is the Euclidean sphere, $\nu$ is the outward unit
normal of  $S_r$ in $\R^3$ and the derivatives are taken with
respect to the Euclidean metric.

\begin{proof}[Proof of Proposition   \ref{examples}]

We first construct metric $g_2$ which contains no compact minimal
surface such that $g_2$ is Schwarzschild  near  infinity. Let
$m>0$ be any positive constant, $\rho_0 > m$, and let $\rho_1
> \rho_0$. Fix a smooth  nonincreasing function $h$ on $(0,
\infty)$ such that $h(r)=0$ for any $r \in (0, \rho_0)$,
$h(r)=-\frac{m}{2}$ for any $r\in (\rho_1, \infty)$. Let $c_0 = 1+
    \frac{m}{2\rho_1}.$
Define a   function $u$ on $\R^3$ as:

\begin{equation}\label{metricg2}
u(x)=c_0 +
 \int_{\rho_1}^{r}\frac{h(\tau)}{\tau^2}d\tau.
\end{equation}
for $|x|=r$. If $|x|=r\le \rho_0$, then
$$u(x)=c_0-\int_{\rho_0}^{\rho_1}\frac{h(\tau)}{\tau^2}d\tau>0
$$ because $h\le
0$. If $|x|=r\in (\rho_0,\rho_1)$, then
$$
u(x)=1+\frac
m{2\rho_1}-\int_{r}^{\rho_1}\frac{h(\tau)}{\tau^2}d\tau>0.
$$
If $|x|=r\ge \rho_1$, then $u(x)=1+\frac{m}{2r}$. Hence $u$ is
smooth and positive. Moreover,
$$
\Delta u=\frac{\p^2 u }{\p r^2}+\frac 2r\frac{\p u}{\p
r}=\frac{h'}{r^2}\le 0.
$$
Let $g_2=u^4g_0$, where $g_0$ is the Euclidean metric. Then
$(\R^3,g_2)$ is an AF manifold which is Schwarzchild near
infinity.

To prove that $(\R^3,g_2)$ contains no compact minimal surfaces,
it is sufficient to show that each Euclidean sphere with center at
the origin has positive mean curvature in $g_2$. Let $  H$ be the
mean curvature of the Euclidean sphere $S_r=\{|x|=r\}$ with
respect to $g_2$. Then

\begin{equation}\label{meancurvature1}
   H=\frac{1}{u^2}\lf(\frac 2r +\frac 4u \frac{\p u}{\p r}\ri)=\frac{2}{ru^3}
   \lf(u+2r\frac{\p u}{\p r}\ri).
\end{equation}
If $r\leq \rho_0$, then  $u'=\frac{\p u}{\p r}=0$. Since  that
$u>0$ is positive,  we have $H>0$. If $r \ge  \rho_0 $, then

 \begin{equation}\label{meancurvature2}
\begin{split}
  2ru' + u &=\frac{2h(r)}{r}+1+\frac{m}{2\rho_1}
  +\int_{\rho_1}^{r}\frac{h(\tau)}{\tau^2}d\tau \\
 &= \frac{2h(r)}{r}+1+\frac{m}{2\rho_1}+h(\xi)\lf(\frac1\rho_1- \frac
 1r\ri)\\
 &\ge \frac{2h(r)}{r}+1+\frac{m}{2\rho_1}+h(r)\lf(\frac1\rho_1 -\frac 1r\ri)
 \\
 &=1+\frac{m}{2\rho_1}+h(r)\lf(\frac1\rho_1 +\frac 1r\ri)\\
 &=1+\frac{m}{2\rho_1}-\frac m2\lf(\frac1\rho_1 +\frac 1r\ri)\\
 &\ge 1-\frac{m}{2\rho_0}\\
 &>0
 \end{split}
\end{equation}
for some $\xi$ between $r$ and $\rho_1$, here   we have used the
fact that $h$ is nonincreasing, $h\ge -\frac m2$ and $r\ge
\rho_0>m$.  Hence we also have $H>0$.

In \cite{M1}, Miao has constructed a scalar flat AF metric $g_1$
on $\R^3$ which contains a stable minimal sphere and is
Schwarzschild at infinity. By rescaling, we see that the
proposition is true. However, because of later application,  we
will construct $g_1$ directly using similar method as in the
construction of $g_2$.

Basically, we   glue the standard sphere to the Schwarzschild
manifold. Under the stereographic projection, the metric on the
sphere minus a point is:
\begin{equation}\label{spherical}
   ds^2_{\mathbb  S^2 }=\frac1{(1+\frac14 \rho^2)^2}(d\rho^2+\rho^2 d\sigma^2)
\end{equation}
where $d\rho^2+\rho^2 d\sigma^2$ is the standard Euclidean metric.
On the other hand the Schwarzschild metric is given by

\begin{equation}\label{Schwarzschild1}
ds^2_{\text{Sch}} = (1+\frac{m}{2r})^4 (d r^2 + r^2 d\sigma^2),
\end{equation}
defined on $\R^3$ minus the origin, where $m>0$ is a constant. We
need to rescale the metric so that the compact minimal surface
$\{r=1/2m\}$ is near $\rho=\infty$  in the metric
(\ref{spherical}). Namely, let $r=\e \rho$, then
\begin{equation}\label{Schwarzschild2}
ds^2_{\text{Sch}} = \e^2(1+\frac{m}{2\e\rho})^4 (d \rho^2 + \rho^2
d\sigma^2),
\end{equation}
For $\rho_0>0$, define
\begin{equation}\label{k}
k(\rho)=\left\{
\begin{array}{ll}
    - \frac \rho4(1+\frac14\rho^2)^{-\frac32}, & \rho \leq \rho_0; \\
    -\frac{m}{2\sqrt\e}\rho^{-2}, & \rho \geq 2 \rho_0. \\
\end{array}\right.
\end{equation}
We want to find $\e>0$ and $\rho_0$ so that $k$ can be defined to
be nonincreasing. Let $\e=m^2/64$, $\rho_0=1/(4m)$. Then
$$
k(\rho_0)=-\frac{32m^2}{(1+64m^2)^\frac 32}.
$$
$$
k(2\rho_0)=-64m^2.
$$
Hence there exists $m_0>0$ such that $k(\rho_0)>k(2\rho_0)$ for
all $0<m<m_0$. For such $m$ we can define $k$ satisfying
(\ref{k}), and is smooth and nonincreasing. Next define
\begin{equation}\label{glue}
  u_m(x)=b_0+\int_0^\rho k(\tau)d\tau
\end{equation}
if $|x|=\rho$, where $b_0$ is chosen such that
$u_m(x)=\sqrt{\e}(1+m/(2\e\rho))$ for $\rho\ge2\rho_0$. More
precisely,
\begin{equation}\label{c0}
\begin{split}
   b_0&=\frac{m}{4\sqrt\e
   \rho_0}+\sqrt\e-\int_0^{2\rho_0}k(\tau)d\tau\\
   &=m\lf(\frac{65}8-\frac{8}{\sqrt{1+64m^2}}\ri)
   +1-\int_{\rho_0}^{2\rho_0}k(\tau)d\tau
   \end{split}
\end{equation}
Note that
\begin{equation}\label{h}
    0\ge\int_{\rho_0}^{2\rho_0}k(\tau)d\tau\ge -16m
\end{equation}
 Hence
\begin{equation}\label{um}
   u_m(x)=\left\{
                  \begin{array}{ll}
                    m\lf(\frac{65}8-\frac{8}{\sqrt{1+64m^2}}\ri)
   -\int_{\rho_0}^{2\rho_0}k(\tau)d\tau+(1+\frac14\rho^2)^{-\frac12},
& \hbox{$\rho\le\rho_0$;} \\
                   m\lf(\frac{65}8-\frac{8}{\sqrt{1+64m^2}}\ri)
   +1-\int_{\rho_0}^{\rho}k(\tau)d\tau, & \hbox{$\rho_0\le 2\rho_0$;} \\
                     \sqrt{\e}(1+\frac m{2\e\rho}), &
\hbox{$\rho\ge2\rho_0$,}
                  \end{array}
                \right.
\end{equation}
for $|x|=\rho$. Hence choosing a smaller $m_0>0$, we have
$u_m(x)>0$ for all $x$ if $0<m<m_0$. As before, since $k<0$ and
$k'\le 0$, we have $\Delta u_m\le0$. The metric
$ds_m^2=u_m^4(d\rho^2+\rho^2d\sigma^2)$ is an AF metric with
nonnegative scalar curvature such that on $\rho\ge2\rho_0$, the
metric is
$$
\e^2(1+\frac{m}{2\e\rho})^4(d\rho^2+\rho^2d\sigma^2)
=(1+\frac{m}{2r})^4(dr+r^2d\sigma^2)
$$
which is Schwarzschild. Moreover, there is a minimal sphere at
$$
\rho=\frac m{2\e}=\frac{32}m>2\rho_0.
$$
Hence $\R^3$ with the metric $ds_m^2$ has a horizon and is
Schwarzschild at infinity.
\end{proof}

From  the proposition, in order to find a sufficient condition for
the existence of compact minimal surfaces, we need to know
information in the interior of the domain. This motivates us to
introduce the following quantity using Hawking mass of some
compact surfaces inside a domain.

   Let $E$ be an
open set in a Riemannian manifold with compact $C^1$ boundary,
  then one can define the mean curvature $H$ of $\p E$ in the weak
  sense, see \cite{HI}. Recall that the Hawking mass $m_H(\p E)$ of $\p  E$
  is defined as:

  \begin{equation}\label{hawking}
 m_H(\p E)=\sqrt{\frac{|\p
  E| }{ 16\pi }}\lf(1-\frac{1}{16\pi}\int_{\p E}H^2\ri).
\end{equation}
where $|\partial E|$ is the area of $\partial E$. In our
convention, the standard unit sphere in $\R^3$ has mean curvature
2 with respect to the unit outward normal.

Let us recall the idea of minimizing hull introduced in \cite{HI}.
Let $(\Omega,g)$ be a Riemannian manifold.

\begin{defi} Let $E$ be a  set in $\Omega$ with locally finite
perimeter. $E$ is said to be a minimizing hull in $\Omega$ if
$|\p^*E\cap K|\le |\p^*F\cap K|$ for any   set $F$ with locally
finite perimeter such that $F\supset E$ and $F\setminus
E\subset\subset \Omega$ and for any compact set $K\subset \Omega$.
Here $\p^*E$ and $\p^*F$ are   the reduced boundaries of $E$ and
$F$ respectively. $E$ is said to be strictly minimizing hull if
equality (for all $K$) implies $E\cap \Omega=F\cap \Omega$ a.e.
\end{defi}
Suppose $E$ is an open set of $\Omega$ such that there is a
strictly minimizing hull in $\Omega$ containing $E$, then define
$E'$ to be the intersection of all strictly minimizing hulls
containing $E$. Note that $E'$ consists of the Lebesque points of
the intersection by definition. $E'$ is called the strictly
minimizing hull of $E$. Let $E\subset\subset \Omega$ be an open
set. Suppose $\Omega$ is compact with smooth boundary which has
positive mean curvature with respect to the outward normal, then
$E'$ exists and $E'\subset \subset \Omega$.

Let $\Omega_1 \subset\subset \Omega_2 \subset \Omega$ such that
$\Omega_1$ and $\Omega_2$ have smooth boundaries. We need the
following lemma from \cite{MY}.

\begin{lemm}\label{areabound}[Meeks-Yau] With the above
notations and let $d$ be the distance between $\Omega_1$ and
$\p\Omega_2$. Let $\iota$ be the infinmum of the injectivity
radius of points in $\{x|\ d(x,\p\Omega_2)>\frac d4\}$. Let $K>0$
be the upper bound of the  curvature of $\Omega_2$. Suppose $N$ is
a minimal surface and $x\in N$ with $d(x,\p\Omega_2)=\frac d2$, so
that $d(x,\p N)\ge \frac d2$, then
\begin{equation}\label{area}
   |N\cap B_x(r)|\ge CK^{-2}\int_0^r \tau^{-1}(\sin
   K\tau)^2d\tau
\end{equation}
where $r=\min\{\frac d2,\iota\}$. Here $C$ is a positive absolute
constant.
\end{lemm}
For such $\Omega_1$, $\Omega_2$, let
\begin{equation}\label{constant}
\alpha^2_{\Omega_1;\Omega_2}=\min\lf\{\frac{CK^{-2}\int_0^r
\tau^{-1}(\sin
   K\tau)^2d\tau}{|\p\Omega_1|},1\ri\}
\end{equation}

Let $\cF_{\Omega_2}$ be the family of precompact connected
minimizing hulls with $C^2$ boundary in $\Omega_2$. Define
\begin{equation}\label{relativelocalmass}
    m(\Omega_1;\Omega_2)=\sup_{E\in \cF_{\Omega_2}, E\subset\Omega_1}m_H(E).
\end{equation}
  Define
\begin{equation}\label{newmass}
    m(\Omega)=\sup  \alpha_{\Omega_1;\Omega_2}m(\Omega_1;\Omega_2)
\end{equation}
where the supremum is taken over all $\Omega_1 \subset\subset
\Omega_2 \subset \Omega$ with smooth boundaries. Here and below
$\Omega_i$ is always assumed to be nonempty.

  In general, the Hawking mass
of a compact surface may be  negative. However, one can   prove
that $m(\Omega)\ge0$.

\begin{theo}\label{nonnegative}
Let $(\Omega,g)$ be a compact manifold with smooth boundary. Then
$\cF_{\Omega_2}\neq \emptyset$ for any $\Omega_2\subset\subset
\Omega$ and $m(\Omega)\ge0$. If $\Omega$ can be embedded in
$\R^3$, then $m(\Omega)=0$. On the other hand, suppose $\Omega$
has nonnegative scalar curvature and $m(\Omega)=0$, then $\Omega$
is locally flat. In particular, if $\Omega$ is simply connected,
then $\Omega$ is a domain in $\R^3$.
\end{theo}

\begin{proof} Since $\Omega$ has smooth
 boundary, by taking a collar
of $\p\Omega$, we can embed $\Omega$ in a compact manifold
$\Omega_3$ with smooth boundary. Taking a double of $\Omega_3$, we
may assume that $\Omega$ is isometrically embedded in a compact
manifold $\Omega_4$ without boundary. Take a small geodesic disk
in $\Omega_4\setminus \ol\Omega$ and glue it to the exterior of
some compact set of $\R^3$, we may assume that $\Omega$ is
embedded in a complete noncompact manifold $M$ with only one end
and such that near infinity of $M$ is isometric to the exterior of
a compact set in $\R^3$.

Suppose $\Omega_2\subset\Omega$. Take a  point $x_0\in \Omega_2$
and let $3r>0$ be such that $B_{x_0}(3r)\subset \Omega_2$,  $3r$
is less than the injectivity radius of $x_0$ and $\p
B_{x_0}(\rho)$ has positive mean curvature for all $0<\rho<3r$. We
claim that there exists $0<\e_0<r$ such that $B_{x_0}(\rho)$ is a
strictly minimizing hull in $\Omega_2$ for all $0<\rho<\e_0$. Let
$0<\rho<r$ and let $F$ be the strictly minimizing hull of
$B_{x_0}(\rho)$ in $M$. $F$ exists because $M$ is Euclidean near
infinity. Moreover, there is $R>0$ independent of $\rho$ such that
$F\subset B_p(R)$ where $p\in M$ is a fixed point. Suppose $F\cap
( M\setminus B_{x_0}(2r))\neq\emptyset$. Then there is a point
$x\in\p F$ and $x\in B_p(2R)\setminus B_{x_0}(2r)$. Moreover,
$B_x(r)\cap   B_{x_0}(r)  =\emptyset$. Since $\p F\setminus
\ol{B_{x_0}}(r)$ is minimal surface, see \cite{HI}, by Lemma 1.1
we have $|\p F|\ge c>0$ for some constant $c>0$ independent of
$\rho$. Now choose $r>\e_0>0$ such that $|\p B_{x_0}(\rho)|<c$ for
all $0<\rho<\e_0$. Then for $0<\rho<\e_0$, the strictly minimizing
hull of $B_{x_0}(\rho)$ must be a subset of $B_{x_0}(3r)$.
However, since $\p B_{x_0}(s)$ has positive mean curvature for all
$0<s<3r$, we conclude that $F=\p B_{x_0}(\rho)$. This proves the
claim. In particular, $B_{x_0}(\rho)\in \cF_{\Omega_2}$ for all
$0<\rho<\e_0$ and so $\cF_{\Omega_2}\neq\emptyset$.

Let $H_\rho $ be the mean curvature of $\p B_{x_0}(\rho)$, then it
is easy to see  that $H_\rho =\frac2\rho +O(1)$ as $\rho\to0$ and
$$
\lim_{\rho\to0}\frac{|\p B_{x_0}(\rho)|}{4\pi\rho^2}=1. $$ Hence
we have $\lim_{\rho\to0}m_H(\p B_{x_0}(\rho))=0$. From this it is
easy to see that $m(\Omega_1;\Omega_2)\ge0$ for all
$\Omega_1\subset\subset\Omega_2\subset\Omega$, and so
$m(\Omega)\geq 0$.

It is well know that $m_H(\Sigma)\le0$ for any compact surface in
$\R^3$, see \cite{W}. Hence if $\Omega$ is a domain in $\R^3$,
then $m(\Omega)=0$.

 To prove the last two assertions,  we first observe that by the
 definition of $m(\Omega)$, $m(\Omega')\le m(\Omega)$
 if $\Omega'\subset \Omega$. Now suppose $m(\Omega)=0$, then for any $x_0\in \Omega$, we have
$m(B_{x_0}(\rho))=0$ for all $\rho$ small enough. We want to prove
$\Omega$ is flat near $x_0$ using the idea in \cite{HI}. Suppose
it is not flat near $x_0$, let $M$ be as above. Let $\rho_0>0$ be
such that (i) $\rho_0<r$ where $3r$ is the injectivity radius of
$x_0$; (ii) $\p B_{x_0}(\rho)$ has positive mean curvature for all
$0<\rho<3r$; (iii)  for all $0<\rho<\rho_0$, $m(B_{x_0}(\rho))=0$
and $B_{x_0}(\rho)$ is a strictly minimizing hull in $M$. Moreover
by the proof above we may assume that for any open set $E\subset
B_{x_0}(\rho_0)$ with $C^2$ boundary, then the strictly minimizing
hull $E'$ of $E$ in $M$ satisfies $E'\subset B_{x_0}(\rho_0)$.

Let $\rho_0>\rho>0$ and let $E_t^\rho$, $0\le t<\infty$ be the
weak solution of inverse mean curvature flow   with initial
condition $B_{x_0}(\rho)$ given by a locally Lipschitz proper
function $u^\rho$. Such a solution exists by \cite{HI}. Moreover,
for $x\in B_{x_0 }(3r)$

\begin{equation}\label{smallball}
|\nabla u^\rho(x)|\le C_1d^{-1}(x)
\end{equation}
for some constant $C_1$ independent of $\rho$ and $d(x)$ is
the distance between $x$ and $x_0$. By
\cite[p.421-422]{HI}, we can find $c_i\to\infty$,
$\rho_i\to 0$ such that $N_{t+c_i}^{\rho_i}$ converges for
a.e. $t\in (-\infty,\infty)$ in $C^1$ to $N_t$ which is a
solution of inverse mean curvature flow in
$M\setminus\{x_0\}$
  given by a locally Lipschitz proper
 function $u$. Here $N_{t+c_i}^{\rho_i}=\p
E_{t+c_i}^{\rho_i}=\p\{u^{\rho_i}<t+c_i\}$ and $N_t=\p
E_t=\p\{u<t\}$. Moreover, $N_t$ is nearly equal to $\p
B_{x_0}(e^{\frac t2})$ as $t\to-\infty$. Since
$N_{t+c_i}^{\rho_i}$ is connected for all $i$,
 $N_{t}$ is also connected.

Let $b>-\infty$ be such that $E_{t+c_i}^{\rho_i}\subset\subset
B_{x_0}(\rho_0)$ for all $t\le b$ and for all $i$. Then as in
\cite[p.427]{HI}, one can show that $m_H(N_{t})>0$ for all $t<b$.
Let us outline the idea. Let  $t<b$ be such that
$N_{t+c_i}^{\rho_i}$
  converges to
 $N_{t}=\p E_{t}$ in $C^1$. Observe that $m_H(N_{t+c_i}^{\rho_i})\geq
m_H (\p
 B_{x_0}(\rho_i))\geq-\delta_0$, here $\delta_0$ is a constant
 independent of $i$
  by the monotonicity formula   \cite[Theorem 5.8]{HI}
 because $\Omega$ has nonnegative scalar curvature.  Since  $|N_{t+c_i}^{\rho_i}|$
 also converges to $|N_{t }|$,  we have
\begin{equation}\label{meancurvature}
   \int_{N_{t+c_i}^{\rho_i}}H^2\le C_2
\end{equation}
for some $C_2$ independent of $i$.
  Since $N_{t+c_i}^{\rho_i}$ converges to $N_{t}$ in
$C^1$, the topology of $N_{t_0+c_i}^{\rho_i}$ is the same
as that $N_{t_0}$ provided $i$ is large enough. By
(\ref{meancurvature}),
  \cite[Lemma 5.5]{HI}   and its proof,  we have
\begin{equation}\label{2ndfundamental}
 \int_{N_{t+c_i}^{\rho_i}}|A|^2\le C_3
\end{equation}
for some constant $C_3$ for all $i$ large enough, where $A$
is the second fundamental form of $N_{t+c_i}^{\rho_i}$. By
lower semicontinuity
\begin{equation}\label{2ndfundamental2}
\int_{N_t}|A|^2<\infty.
\end{equation}
Hence the monotonicity formula \cite[Theorem 5.8]{HI} for
$N_t$ holds. Note that by upper semicontinuity
$m_H(N_t)\ge0$ for such $t$. Hence if for such a $t$,
$m_H(N_t)=0$, then one can prove that $m_H(N_s)=0$ for a.e.
$s<t$. Since there is no compact minimal surfaces in
$B_{x_0}(3r)$, one can argument as in \cite{HI} that
$\Omega$ is flat near $x_0$. This is a contradiction.

Now choose $t_0<b$ such that $N_{t_0+c_i}^{\rho_i}$
converge in $C^1$ to $N_{t_0}$. Note that even though
$m_H(N_{t_0})>0$  and $E_{t_0}$ is a
 minimizing hull by \cite{HI}, $N_{t_0}$ may not have $C^2$
 boundary. In order to show that $m(B_{x_0}(\rho))>0$ for $\rho$ small, we need some
 approximation.

As before, we have
$$
\int_{N_{t_0}}|A|^2<\infty.
$$
By (5.19) in
 \cite{HI}, we can find an open set  $F\subset B_{x_0}(\rho_0) $ with smooth
boundary such
 that   $m_H(\p F)>0$, $\p F$ is
 connected.
Let
  $F'$ be the
strictly minimizing hull of $F$. Then
  $F'\subset B_{x_0}(\rho_0)$ and
\begin{equation}
   \begin{split}
m_H(\p F')&=\sqrt{\frac{|\p
  F'| }{ 16\pi }}\lf(1-\frac{1}{16\pi}\int_{\p F'}H^2\ri)\\
  &\ge
  \sqrt{\frac{|\p
  F'| }{ |\p F| }}m_H(\p F)\\
  &>0
\end{split}
\end{equation}
because the mean curvature of $\p F'$ is zero on $\p F'\setminus\p
F$ and the mean curvature of $\p F'$ is equal to the mean
curvature of $\p F$ a.e. on $\p F'\cap\p F$.

   Since $B_{x_0}(3r)$ is foliated by
  geodesic sphere with positive mean
  curvature,   $\p F$
     cannot be a minimal surface. Since $F$ is connected,
     $F'$ is connected. Since $\p F'$ is $C^{1,1}$,
     by Lemma 5.6 in \cite{HI}  we can find
  $G\subset B_{x_0}(\rho_0)$ with smooth boundary such that $G$ is a
  strictly minimizing hull in $M$
  and $m_H(G)>0$.  $G$ can be chosen to be
  connected because $F'$ is connected. Therefore $G\in\cF_{B_{x_0}(\rho_0)}$, we have
  $m( B_{ x_0 }(\rho_0))>0$. This is a contradiction because $m(B_{x_0}(\rho_0))=0$. Hence
$\Omega$ is locally flat. If $\Omega$ is simply connected, then
$\Omega$ is a domain in $\R^3$, see \cite[p.42-44]{Sp}.

\end{proof}

\section  { Sufficient conditions on the existence of compact minimal
surfaces}

Let $(\Omega  ,g)$ be a compact orientable three manifold with
smooth boundary which consists of finitely many components and
with nonnegative scalar curvature.  In this section, we always
assume that the boundary of $\Omega$ has positive Gauss curvature
and positive mean curvature with respect to the outward normal.
Hence each boundary component is diffeomorphic to a sphere. We
want to obtain a sufficient condition on the existence of {\it
outermost horizon} which is defined as in \cite{B} as follows:
  Let $\Sigma_1,\dots,\Sigma_k$ be the
components of $\p\Omega$.   For a fixed $i$, by adding points to
each component  except $\Sigma_i$, we obtain a topological space
$\Omega_i$. Let $\mathcal S_i$ be the collection of surfaces which
are smooth boundaries of precompact open sets in $\Omega_i$
containing those points. The boundary $\p E$ of a set  $E$ with
$\p E\in \mathcal S_i$ is called outer minimizing if $|\p E|\le
|\p F|$ for all $F\supset E$. A  {\it horizon} (relative to
$\Sigma_i$) is a surface $\p E\in \mathcal S_i$  with zero mean
curvature. A horizon is said to be {\it outermost} if it is not
enclosed by another horizon, see \cite[p.185]{B}. Suppose $\p E$
is a horizon. Let $E'$ be its strictly minimizing hull in
$\Omega_i$, which exists because the mean curvature of $\Sigma_i$
is positive. Then $\p E'$ is also a  horizon.

 Let $(\Omega,g)$ and $\Sigma_i$ as before.  By Weyl embedding theorem,
  see \cite{N}, each $\Sigma_i$ can
be isometrically embedded in $\R^3$ and the embedding is unique up
to an isometry of $\R^3$. Let $H$ be the mean curvature of
$\Sigma_i$ with respect to the outward normal and let $H_0$ be the
mean curvature of $\Sigma_i$ when embedded in $\R^3$. The
Brown-York mass \cite{BY1,BY2} of $\Sigma_i$ is defined to be
\begin{equation}\label{BYmass}
   m_{BY}( \Sigma_i )=\frac{1}{8\pi}\int_{ \Sigma_i}(H_0-H)d\sigma
\end{equation}
where $d\sigma$ is the volume element on $\Sigma_i$ induced by
$g$.

In this section we will give some  lower bounds of the Brown-York
mass in terms of the Hawking mass of subsets of $\Omega$. We will
also give a sufficient condition on  of existence of horizons in
$\Omega$ using the Brown-York mass and $m(\Omega)$ defined in the
previous section.  The first lower bound for the Brown-York mass
is the following:

\begin{theo}\label{lowerboundBYmass}  Let $(\Omega,g)$ be a compact
three manifold with connected smooth boundary and with nonnegative
scalar curvature. Assume that $\Omega$ is simply connected and
suppose $\p\Omega$ has positive Gauss curvature and positive mean
curvature with respect to the outward normal.  Then
\begin{equation}
    m_{BY}(\p \Omega)\ge   m_H(\p E).
\end{equation}
for any   connected minimizing hull $E$ in $\Omega$ where
$E\subset\subset \Omega$   with $C^{1,1}$ boundary.  Moreover,
equality holds for some minimizing hull $E$ with the above
properties if and if $\Omega$ is a standard ball in $\R^3$ and $E$
is a standard ball in $\Omega$. In particular,
$m_{BY}(\p\Omega)\ge m_H(\p\Omega)$ and equality holds if and only
if $\Omega$ is a standard ball in $\R^3$.
\end{theo}
\begin{proof} Let $E$ be a minimizing hull in $\Omega$ with
$C^{1,1}$ boundary such that $E\subset\subset\Omega$. We want to
prove that $m_{BY}(\p\Omega)\ge m_H(\p E)$. Since
$m_{BY}(\p\Omega)\ge0$,  it is sufficient to prove the case when $
m_H(\p E)>0$.

 Isometrically embed
$\Sigma=\p\Omega$ in $\R^3$. As in \cite{ST}, one can glue
$\Omega$ to the exterior of $\p\Omega$ in $\R^3$ to form a
manifold $(M,h)$ such that the metric $h$ satisfies:

\begin{enumerate}
    \item [(i)] $h|_\Omega=g$;
    \item [(ii)] $h_{M\setminus\Omega}$ is smooth up the
boundary;
    \item [(iii)] the scalar curvature of $h$ in $M\setminus \Omega$
is zero;
    \item [(iv)] $h$ is Lipshitz near $\p\Omega$;
    \item [(v)] the mean
curvatures of $\p\Omega$ with respect to the outward normal are
the same for metrics inside and outside $\Omega$;
    \item [(vi)] $h$ is asymptotically flat;
    \item [(vii)] if $\Sigma_r$ is
the surface consisting points in the exterior to $\Sigma$ in
$\R^3$ and with Euclidean distance $r$ from $\Sigma$, then
$\Sigma_r$ has positive mean curvature with respect to $h$, which
is bounded from below by $C/(1+r)$ for some $C>0$.
\end{enumerate}

 Given $\e>0$, by (i)-(vi) and   \cite{M},   there is a smooth
 AF metric $h_\e$ on $M$ with nonnegative scalar curvature such that
 \begin{equation}\label{metric}
    (1-\e)h_\e(v,v)\le h(v,v)\le (1+\e)h_\e(v,v)
\end{equation}
for all tangent vector $v$ on $M$, and such that
\begin{equation}\label{mass}
   \lim_{\e\to0}m_{ADM}(h_\e)=m_{ADM}(h)
\end{equation}
where $m_{ADM}(h_\e)$ and $m_{ADM}(h)$ are the ADM masses of
$(M,h_\e)$ and $(M,h)$ respectively.

Let $\theta>0$ be given. We can find a connected open set
$F\supset E$ with smooth boundary such that $F
\subset\subset\Omega$,
\begin{equation}\label{approximation}
   |\p E|-\theta\le |\p F|\le |\p E|+\theta; \ m_H(\p E)\le m_H(\p
    F)+\theta; \text{and}\ m_H(\p F)>0.
\end{equation}

For any $\e>0$, let $F'$ be the strictly minimizing hull of $F$ in
$(M,h_\e)$. $F'$ exists, precompact in $M$  and has $C^{1,1}$
boundary. Since $F$ is connected, $F'$ is connected. Moreover $M$
is simply connected because $\Omega$ is simply connected and
$\p\Omega$ is homeomorphic to the sphere. Hence $\p F'$ is
connected.

Choose $\delta>0$ be small enough, such that the boundary of
$\Omega_t$ for all $0<t\le \delta$ has positive mean curvature in
$g$, where $\Omega_t=\{x\in \Omega| \ d_g(x,\p\Omega)>t\}$. Here
$d_g$ is the distance with respect to $g$. There is $0<t_0<\delta$
such that $F\subset\subset \Omega_{t_0}$. Suppose $F'\setminus
\ol{\Omega_{t_0}}\neq\emptyset$.  Since $M\setminus \Omega_\delta$
is foliated by compact surfaces of positive mean curvature with
respect to $h$, if we denote the unit outward normal vector in $h$
of the surfaces by $\nu$, we have
\begin{equation}\label{comparison1}
\begin{split}
0&<\int_{F'\setminus\Omega_{t_0}} {div}_h\nu dV_h\\
&=\int_{\p F'\cap (M\setminus \ol{\Omega_{t_0}})}h(\nu,\mu)-|\ol
F'\cap \p\Omega_{t_0}|_g
\end{split}
\end{equation}
where $\mu$ is the unit outward normal of $\p F'$ in $h$. Here $
{div}_h$ is the divergence with respect to $h$.  Hence we have
\begin{equation}\label{comparison2}
 |\p F'\cap (M\setminus \Omega_{t_0})|_h> |\ol F'\cap
\p\Omega_{t_0}|_g
\end{equation}
Now $F'\cap \Omega_{t_0}\supset E$ is a set of finite perimeter
and $E$ is a minimizing hull in $(\Omega,g)$, we have
\begin{equation}\label{comparison3}
  |\p (F'\cap \Omega_{t_0})|_g\ge |\p E|_g.
\end{equation}
  From (\ref{comparison2}) and (\ref{comparison3}) we
conclude that
\begin{equation}\label{comparison4}
\begin{split}
  |\p F'|_h&=|\p F'\cap (M\setminus \Omega_{t_0})|_h
  +|\p F'\cap \Omega_{t_0}|_h\\
  &> |\ol F'\cap
\p\Omega_{t_0}|_g+|\p F'\cap \Omega_{t_0}|_g\\
&=|\p (F'\cap \Omega_{t_0})|_g\ge |\p E|_g.
\end{split}
\end{equation}
Moreover the mean curvature of $\p F'$ in $h_\e$ is zero on $\p
F'\setminus \p F$ and the mean curvature of $\p F'$ and $\p F$ are
equal a.e. in $\p F'\cap \p F$. Using  the inverse mean curvature
flow with initial data $F'$, we obtain by the Penrose inequality
in \cite{HI}  we have
\begin{equation}\label{penrose}
\begin{split}
    m_{ADM}(h_\e)&\ge \sqrt{\frac{|\p
  F'|_{h_\e} }{ 16\pi }}\lf(1-\frac{1}{16\pi}\int_{\p F'}H^2d\sigma_{h_\e}\ri)
\\
  & \ge \sqrt{\frac{|\p
  F'|_{h_\e} }{|\p F|_{g}}}\cdot\sqrt{\frac{|\p
  F |_{g} } {16\pi}}\lf(1-\frac{1}{16\pi}\int_{\p F}H^2d\sigma_{h_\e}\ri) \\
  &=\sqrt{\frac{|\p F'|_{h_\e} }{|\p F|_{g}}}m_H(\p F)\\
&\ge \sqrt{\frac{|\p F'|_{h_\e} }{|\p F|_{g}}}(m_H(\p E)-\theta)\\
&\ge \sqrt{(1-\e)\frac{|\p F'|_{h} }{|\p F|_{g}}}m_H(\p E)-\theta
\sqrt{\frac{|\p F|_{h_\e} }{|\p F|_{g}}}\\
&\ge \sqrt{(1-\e)\frac{|\p E|_{g} }{|\p F|_{g}}}m_H(\p
E)-\theta(1+\e)^\frac12\\
&\ge \sqrt{(1-\e)\frac{|\p E|_{g} }{|\p E|_{g}+\theta}}m_H(\p
E)-\theta(1+\e)^\frac12
\end{split}
\end{equation}
Here, we have use the assumption $m_H (\p E)\ge 0$ in the fifth
inequality above, and  $d\sigma_{h_{\e}}$ and $d\sigma_g$ are the
area elements on surfaces induced by $h_\e$ and $g$ respectively
and we have used (\ref{approximation}),  and (\ref{comparison4}).
Let $\e\to0$, and then let $\theta\to0$, we have $m_{ADM}(h)\ge
m_H(\p E)$. By \cite{ST}, we have
\begin{equation}\label{comparison5} m_{BY}(\p\Omega)\ge
m_{ADM}(h)\ge m_H(\p E)
\end{equation}
If $m_{BY}(\p\Omega)= m_H(\p E)$ for some $E$ with the properties
in the theorem, then $m_{BY}(\p\Omega)= m_{ADM}(h)$ by
(\ref{comparison5}). Hence $\Omega$ is a domain in $\R^3$ by
\cite{ST} and $m_H(\p\Omega)=m_{BY}(\p\Omega)=0$. By    \cite{W},
$\Omega$ is a standard ball in $\R^3$. Since $m_H(\p E)=0$, $E$ is
also a standard ball in $\R^3$.

 Conversely, if
$\Omega$ is a domain in $\R^3$, since each sphere has zero Hawking
mass, and a ball which is small enough is a minimizing hull, we
have $m_{BY}(\Omega)= m_H(\p E)=0$ for some small ball.

To prove the last statement, we observe that $ \Omega_t$ is a
minimizing hull in  $\Omega$ if $t>0$ is small. Hence
$m_{BY}(\p\Omega)\ge\lim_{t\to0}m_H(\p\Omega_t)=m_H(\p\Omega)$.
In fact, from the proof above,
$m_{BY}(\p\Omega)>\lim_{t\to0}m_H(\p\Omega_t)$ unless $\Omega$ is
a domain in $\R^3$. Hence by a result of Willmore \cite{W},
$m_{BY}(\p\Omega)=m_H(\p\Omega)$ implies that $\p\Omega$ is a
standard sphere so $\Omega$ is a standard ball in $\R^3$.
\end{proof}

In the above theorem, we do not assume that $(\Omega,g)$  contains
no horizons. In order to obtain a sufficient condition that
$(\Omega,g)$ contains a horizon, we need another estimate of the
Brown-York mass. Let $m(\Omega)$ as defined in (\ref{newmass}). We
have the following:

\begin{theo}\label{blackhole}  Let $(\Omega,g)$ be a compact
manifold with connected smooth boundary and with nonnegative
scalar curvature. Assume that $\Omega$ is simply connected and
suppose $\p\Omega$ has positive Gauss curvature and positive mean
curvature with respect to the outward normal. Suppose $(\Omega,g)$
has no    horizons. Then  $m(\Omega)\le m_{BY}(\p\Omega)$.
Equality holds if and only if $\Omega$ is a domain in $\R^3$.
\end{theo}
\begin{rema} If $\p\Omega$ has more than one components, then it
is easy to see that $\Omega$ will have an outermost minimizing
horizon with respect to any boundary component because $\p\Omega$
has positive mean curvature. If $\p\Omega$ is connected and
$\Omega$ is not simply connected, then $\Omega$ contains a minimal
  sphere or real projective space by \cite{MSY}.
\end{rema}

\begin{proof}[Proof of Theorem \ref{blackhole}]
 Assume that
$\Omega$ has   no horizons. Suppose $ \Omega$ is not diffeomorphic
to an open 3-ball in $\R^3$, then $\Omega$ must contains a minimal
sphere $S$ or projective plane \cite{MSY}. Since $\Omega$ is
simply connected, $S$ must be a minimal sphere  and is the
boundary of some precompact open set $E$, see  \cite[p.107]{H}. So
$S$ is a horizon. This is a contradiction. Hence $\Omega$ is
diffeomorphic to an open 3-ball in $\R^3$.

 Let $\Omega_1\subset\subset\Omega_2\subset \Omega$ with smooth
 boundaries. Let $E$ be a connected  precompact minimizing hull in
 $\Omega_2$ with $C^2$ boundary such that $E\subset\Omega_1$. We want to prove that
$m_{BY}(\p\Omega)\ge \alpha_{ \Omega_1;\Omega_2}m_H(\p E)$. Again,
we may assume that $m_H(\p E)\ge 0$.

 Let $E'$ be the
 strictly minimizing hull of $E$ in $\Omega$.   $E'$ exists
 because $\p\Omega$ has positive mean curvature. Moreover, $E'$
 has $C^{1,1}$ boundary because $E$ has $C^2$ boundary.
Note that $E'$ and $\p E'$
 are connected because $\Omega$ is diffeomorphic to the unit ball in $\R^3$.
  By
 Theorem \ref{lowerboundBYmass},
 \begin{equation}\label{estimates1}
 \begin{split}
 m_{BY}(\p\Omega)&\ge m_H(\p E')\\
 &\ge \sqrt{\frac{|\p
  E'| }{ 16\pi }}\lf(1-\frac{1}{16\pi}\int_{\p E}H^2\ri)
  \end{split}
\end{equation}
because $\p E'\setminus \p E$ is minimal and the mean curvatures
of $\p E$ and $\p E'$  are equal a.e. on their common part.
Suppose $E'\subset\subset \Omega_2$, then $|\p E'|\ge |\p E|$ and
we have $m_{BY}(\p\Omega)\ge m_H(\p E)$.

Suppose $E'$ is not a precompact set in $\Omega_2$, then $\p
E'\cap( M\setminus \Omega_2)\neq\emptyset$. On the other hand, if
$\p E'\cap\p E=\emptyset$ then  $\p E'$ would be a horizon
contradicting the assumption that $\Omega$ contains no horizons.
Hence $$\p E'\cap \Omega_1\supset\p E'\cap\p E\neq\emptyset$$ and
there is a point $x\in \p E'$ and $d(x,\p\Omega_2)=\frac d2$,
where $d$ is the distance between $\Omega_1$ and $\p\Omega_2$. By
Lemma \ref{areabound}, we have
\begin{equation}\label{estimates2}
   |\p E'|\ge \alpha^2_{ \Omega_1;\Omega_2}|\p\Omega_1|\ge
   \alpha^2_{\Omega_1;\Omega_2}|\p E|
\end{equation}
because $E$ is a minimizing hull in $\Omega_2$. Hence by
(\ref{estimates1}), (\ref{estimates2}) we have
\begin{equation}\label{estimates3}
  m_{BY}(\p \Omega)\ge \alpha_{\Omega_1;\Omega_2}m_H(\p E).
\end{equation}
This completes the proof of the first part. If equality holds,
then from the proof above, for any $\e>0$ we can find a strictly
minimizing hull $E'$ in $\Omega$ with $C^{1,1}$ boundary such that
$m_{BY}(\p\Omega)\le m_H(\p E')+\e$. From the proof of Theorem
\ref{lowerboundBYmass}, one can conclude that $\Omega$ must be a
domain in $\R^3$.
\end{proof}
Let $\Omega$ be as in the theorem. Isometrically embed $\p\Omega$
 in $\R^3$. Let $R$ be the radius of the smallest
circumscribed ball of $\p\Omega$ in $\R^3$. We have the following:

\begin{coro}\label{blackholecor} Let $(\Omega,g)$ be a compact manifold with
nonnegative scalar curvature and with connected boundary which has
positive mean curvature and positive Gauss curvature. Suppose
$\Omega$ is simply connected and suppose $m(\Omega)\ge
m_{BY}(\p\Omega)$, then there is a horizon in $\Omega$ unless
$\Omega$ is a domain in $\R^3$. Hence if $m(\Omega)\ge 2R$, then
$\Omega$ contains a horizon. In particular, if $m(\Omega)\ge
2\text{diam}(\p \Omega)$ then $\Omega$ contains a horizon. Here
$\text{diam}(\p \Omega)$ is the diameter of $\p\Omega$ with
respect to the metric induced by $g$.
\end{coro}
\begin{proof} The first statement is just a restatement of
Theorem \ref{blackhole}. Since
\begin{equation}\label{diam}
m_{BY}(\Omega)\le \frac{1}{8\pi}\int_{\p \Omega}H_0 \le 2R <
2\text{diam}(\p\Omega)
\end{equation}
by the Minkowski integral formula \cite[p.136]{K}, the fact that
$\p\Omega$ has positive Gauss curvature and the Gauss-Bonnet
formula.
\end{proof}

 There are   examples so that $m(\Omega)>m_{BY}(\Omega)$.
 Consider the metrics $(\R^3,ds_m^2= u_m^4(d\rho^2+\rho^2d\sigma^2))$
 defined in the proof of Proposition  \ref{examples},
 where $d\rho^2+\rho^2d\sigma^2$ is the Euclidean metric.
 Since as $m\to0$, $ds_m^2$ converges uniformly on the Euclidean ball
  $B_4=\{x\in \R^3|\ |x|<4\}$ to the standard metric on the unit sphere,
   we see that    for sufficiently small $m>0$, the boundary  of
$B_\tau$   is mean convex for all $0<\tau<2$ with respect to
$ds_m^2$.     $B_1$ is a minimizing hull in $B_2$. Moreover, there
is $\delta>0$ such that the Hawking mass $m_H(\p B_1)$ of $\p B_1$
with   respect to $ds_m^2$ is at least $\delta$, provided $m>0$ is
small enough. Hence for all $m>0$ small enough,
$\alpha_{\Omega_1;\Omega_2}m(B_1;B_2)\ge \delta$.  Now consider
the domain $\Omega=B_\rho$, where $\rho=8/m>2\rho_0$. It is easy
to see that $\p\Omega$ has positive mean curvature and positive
Gauss curvature with diameter $ d\leq Cm$, here $C$ is a universal
constant. Hence
$$
m(\Omega)\ge \alpha_{\Omega_1;\Omega_2}m(B_1;B_2)\ge \delta>2d\ge
m_{BY}(\Omega)
$$
provided $m$ is small enough.

There is another sufficient condition for the existence of
horizons in terms of the Hawking mass of other kind of surfaces.
Suppose $E\subset\subset \Omega$ is a bounded subset in $\Omega$
such that $E$ and connected and the volume of $E$ is equal  $v_0$.
Using the terminology in \cite{CY}, $\p E$ is called a {\it round
sphere} if $\p E$ is smooth and connected and  has least area
among  all subsets of $\Omega$ with locally finite perimeter with
volume $v_0$. The volume of a set of locally finite perimeter is
simply the   measure of the set with respect to the volume element
of $\Omega$.
 If the scalar curvature of
$\Omega$ is nonnegative, and $\p E$ is a round sphere, then the
Hawking mass $m_{H}(\p E)$ of $\p E$ is nonnegative by  \cite{CY}.
We have the following:

\begin{theo}\label{blackhole1}Let $(\Omega,g)$ be a compact
Riemannian three manifold with smooth  and connected boundary and
with nonnegative scalar curvature. Assume that $\Omega$ is simply
connected and suppose that  both the Gauss curvature and mean the
curvature of the boundary( with respect to outward normal) is
positive. If there is a round sphere $\p E$ in $\Omega$ such that

$$m_H (\p E)> m_{BY}(\partial \Omega),$$
then there is a horizon in $\Omega$.
\end{theo}

First, we prove the following lemma:

\begin{lemm}\label{roundsphere}   Let   $(\Omega,g)$ be a compact
Riemannian three manifold with smooth    boundary and with
nonnegative scalar curvature such that $\p\Omega$ has positive
mean curvature. Suppose $E\subset\subset\Omega$ such that $\p E$
is a round sphere. Then either $\Omega$ contains a horizon or $E$
is a minimizing hull in $\Omega$.
\end{lemm}
\begin{proof}
Suppose $E$ is not minimizing hull in $\Omega$. Since $\p\Omega$
has positive mean curvature, there exists a strictly minimizing
hull $E'$ containing $E$ with $E'\subset\subset\Omega$. Since $E$
is not a minimizing hull, $|\partial E'| < |\p E|$. Let $v_0$ be
the volume of $E$, then $vol(E')>v_0$ because $E$ is a round
sphere. Choose $\e_0>0$ small enough so that $\p\Omega_\e$ is
smooth and has positive mean curvature for all $0<\e<2\e_0$, where
$\Omega_\e=\{x\in\Omega|\ d(x,\p\Omega)>\e\}$. Moreover, $\e_0$ is
chosen so that $E'\subset\subset\Omega_{\e_0}$. Let $\Phi$ be the
family of precompact subsets $F$ with locally finite perimeter of
$\Omega_{\e_0}$ with the properties that $vol(F)> v_0$,  and $|\p
F|<\frac12( |\p E'|+|\p E|)$. $\Phi$ is not empty  because $E'\in
\Phi$. Let
\begin{equation}\label{inf1}
   v_1= \inf_{F\in\Phi} vol(F).
\end{equation}
Note that $v_1\ge v_0$. We claim that $v_1>v_0$. In fact, let
 $ F_i\in \Phi$ be
such that $\lim_{i\to\infty}vol( F_i)=v_1$.
 By compactness theorem for sets of locally
finite perimeter, we may assume that the characteristic functions
of $F_i$ converges in $L^1(\Omega)$ to a the characteristic
function of a subset $F$ of $\ol{ \Omega}_{\e_0}$ with locally
finite perimeter such that $vol(F)=v_1$. By lower semicontinuity,
$|\p F |\le \frac12( |\p E'|+|\p E|)<|\p E|$. (For simplicity, we
use $|\p F|$ to denote the area of the reduced boundary of $F$).
Hence $v_1=vol(F)>v_0$ by the definition of round spheres.

Let
\begin{equation}\label{inf2}
   a= \inf_{F\in\Phi} |\p F|.
\end{equation}
Choose a sequence $F_i\in \Phi$ such that $\lim_{i\to\infty}|\p
F_i|=a$. Since $v_1>v_0$ and $|\p F_i |< \frac12( |\p E'|+|\p
E|)$, by the approximation result of sets of locally finite
perimeter, see \cite{Gi}, we may choose $F_i$ such that $\p F_i$
is smooth. Choose $2\e_0>\e_1>\e_0$ such that
$vol(\Omega_{\e_0}\setminus \Omega_{\e_1})\le \frac12(v_1-v_0)$.
Consider the set $\tilde F_i=F_i\cap \Omega_{\e_1}$, then $|\p
\tilde F_i|\le |\p F_i|<\frac12( |\p E'|+|\p E|)$  since $\p
\Omega_\e$ has positive mean curvature for all $0<\e<2\e_0$.
Moreover,
\begin{equation}\label{inf3}
vol(\tilde F_i)\ge vol(F_i)-vol(\Omega_{\e_0}\setminus
\Omega_{\e_1})\ge \frac12(v_1+v_0)>v_0
\end{equation}
by (\ref{inf1}).  Hence $\tilde F_i\in \Phi$ and
$\lim_{i\to\infty}|\p \tilde F_i|=a$.  Passing to a subsequence if
necessary, we may assume that   the characteristic functions of
$\tilde F_i$ converges in $L^1(\Omega)$ to a the characteristic
function of a subset $\tilde F$ of $\ol{ \Omega}_{\e_1}$ with
locally finite perimeter. Moreover, $vol(\tilde F)\ge
\frac12(v_1+v_0)>v_0$ by (\ref{inf3}) and $|\p \tilde F|=a$. Since
$E'\in \Phi$, $a\le |\p E'|<\frac12(|\p E'|+|\p E|)$.  Hence
$\tilde F\in \Phi$.

Let $B_x(r)\subset\subset\Omega_{\e_0}$ with
$vol(B_x(r))<\frac12(v_1-v_0)$. Suppose $F$ is any set of locally
finite perimeter such that the symmetric difference $F
\Delta\tilde F\subset\subset B_x(r)$, then
$F\subset\subset\Omega_{\e_0}$,
$$vol(F)\ge vol(\tilde
F)-vol(B_x(r))\ge \frac12(v_1+v_0)>v_0 $$ by (\ref{inf1}). If $|\p
F|<|\p \tilde F|$ then $F\in \Phi$ because $|\p \tilde
F|<\frac12(|\p E'|+|\p E|)$. This is impossible  by the
construction of $\tilde F$. Hence one can see that $\tilde F $ is
minimizing in any geodesic ball inside $\Omega_{\e_0}$ with radius
small enough and  $\tilde F$ is a compact minimal surface. That is
to say $\Omega$ contains a horizon.
\end{proof}

\begin{proof}[  Proof of Theorem \ref{blackhole1}]  Let $\p E$ be a round
sphere.   $E$ is not a minimizing hull in $\Omega$ by Theorem
\ref{lowerboundBYmass} and the assumption that $m_H(\p
E)>m_{BY}(\p\Omega)$. By Lemma \ref{roundsphere}, we see that
$\Omega$ contains a horizon. This completes the proof of the
theorem.
\end{proof}

\section {On  the Bartnik quasi-local mass}

In this section, we will discuss relation between $m(\Omega)$ and
the quasi-local mass introduced by Bartnik \cite{B1}. We will also
discuss some properties of the Bartnik quasi-local mass.
 We use the definition of the
Bartnik quasi-local mass $m_{B}(\Omega)$  in \cite{B2} of a
compact manifold $( \Omega,g)$ with smooth boundary and with
nonnegative scalar curvature.

Let $ \mathcal{PM}_o $  be the set of complete noncompact AF
manifolds with nonnegative scalar curvature which admit no minimal
two spheres or projective planes. For $M\in \mathcal{PM}_o$,
$\Omega\subset\subset M$ means that $\Omega$ is isometrically
embedded in $M$. Then the Bartnik quasi-local mass $m_B(\Omega)$
of $\Omega$ is defined as:
\begin{equation}\label{bartnikmass}
   m_B(\Omega)=\inf_{M\in \mathcal{PM}_o }\{m_{ADM}(M)|\ \Omega\subset\subset M\}
\end{equation}
If $\Omega\subset\subset M\in \mathcal{PM}_o$,
 then $M$ is said to be an admissible extension of $\Omega$.
We should remark that for $M\in \mathcal{PM}_o$,
  $M$ is topologically $\R^3$, by
\cite{MSY}, see also \cite{G}.

We have the following lower bound of $m_B(\Omega)$.

\begin{theo}\label{bmasslowerbound} Suppose $\Omega$ has an
admissible extension.   Then
\begin{equation}\label{blowerbound1}
    m_B(\Omega)\ge  m(\Omega).
\end{equation}
\end{theo}
\begin{proof} The proof is exactly the same as the proof of
 Theorem \ref{lowerboundBYmass}.
 \end{proof}

It was proved in \cite{HI}
  that if
  $m_B(\Omega)=0$, then $\Omega$ is locally flat. Hence it is a
   domain in $\R^3$ if
$\p\Omega$ is topologically a sphere. If the boundary of $\Omega$
has several components, then one can show that if
$m_B(\Omega)=0$, then $\Omega$ is still a domain in $\R^3$
provided that each component of $\p\Omega$ is topologically a
sphere. This follows from the following lemma.

\begin{lemm}\label{positivity1} Let $(\Omega,g)$ be a compact three
 manifold with
boundary with nonnegative scalar curvature such that $\Omega$ has
an admissible extension. Suppose $m_B(\Omega)=0$. Let
$\Sigma_1,\dots,\Sigma_k$ be the boundary components of $\Omega$.
Then $\Omega$ can be extended to a manifold $\tOmega$ with compact
closure with nonnegative scalar curvature such that $\tOmega$ has
only one boundary component which is equal to $\Sigma_i$ for some
$i$ and $m_B(\tOmega)=0$.
\end{lemm}
\begin{proof} Let $M_l$ be a minimizing sequence of admissible
extensions of $\Omega$ for $m_B(\Omega)$. Then each $M_l$ is
diffeomorphic to $\R^3$. Each   $\Sigma_j$ separates $\R^3$ into
two components.  Hence by taking a subsequence if necessary, we
may assume that $\Omega$ is in the interior of $\Sigma_1$, say,
for all $l$. Let $\tOmega$ be the interior of $\Sigma_1$ in $M_1$.
Note that $\tOmega$ is diffeomorphic to the unit ball in $\R^3$.
Then $\Omega\subset \tOmega$, the scalar curvature of $\tOmega$ is
nonnegative, and $\p\tOmega=\Sigma_1$. We want to prove that
$m_B(\tOmega)=0$.

Let $\tM_l$ be the manifold obtained in the following way:
$\tM_l=M_l$ on  $M_l\setminus$interior of $\Sigma_1$ and replace
$\Omega$ by $\tOmega$ in the interior of $\Sigma_1$, so that in a
neighborhood of $\Sigma_1$, $M_l=\tM_l$.  We want to prove that if
$l$ is large, $\tM_l$ will not contain any stable  minimal
spheres.

 Choose a $\e>0$ be small enough, such that
  $\{x\in \Omega|\ d(x,\Sigma_1)=\e\}$ is also topologically a sphere
  and so that if $U= \{x\in \Omega|\ d(x,\Sigma_1)<\e\}$,
   then $\ol U\cap \Sigma_j=\emptyset$ for all $j\neq1$. Let $N$ be an outermost
   horizon of $\tM_l$.
   If $N\cap
(\tOmega\setminus U)=\emptyset$, then $N\subset M_l$. This is
impossible by the fact that $M_l$ is admissible and the fact that
each component of $N$ is a stable minimal sphere, see \cite{HI,B}.
Hence $N\cap( \tOmega\setminus U)\neq\emptyset$ and there is a
constant $c>0$ independent of $i$, such that $|N|\ge c>0$ by Lemma
\ref{areabound}. On the other hand,
$$
\lim_{l\to\infty}m_{ADM}(\tM_l)= \lim_{l\to\infty}m_{ADM}(M_l)=0.
$$
Hence by the Penrose inequality in \cite{HI, B}, no such $N$
exists if $l$ is large enough.

 Note that $\tM_l$ is homeomorphic
to $\R^3$. Hence $\tM_l$ contains no minimal projective plane. If
$\tM_l $ contains a minimal sphere $S$, then strictly minimizing
hull of the interior of $S$ is an outermost horizon. Hence no such
$S$ exists. We conclude that $\tM_l$ is an admissible extension of
$\tOmega$ and so $m_B(\tOmega)=0$.
\end{proof}
From this, we obtain  the following:
\begin{prop}\label{positivity} Let $(\Omega,g)$ be as in the
lemma.
 Suppose $m_B(\Omega)=0$ and each $\Sigma_i$
  is topologically a sphere, then $\Omega$ is a domain in $\R^3$.
\end{prop}

In case the boundary of
 $\p \Omega$ is not a sphere, it is still unknown  if
 $m_B(\Omega)=0$ will imply $\Omega$ is a domain in $\R^3$. In the
 following, using Theorem \ref{bmasslowerbound} we are going to construct a locally
 flat compact manifold which either has no admissible extension or
 its Bartnik mass is positive.

\begin{exam}\label{cylinder} Consider the cylinder in $\Bbb{R}^3$:
$$C=\{(x,y,z)\in \Bbb{R}^3 |x^2 +y^2 \leq r^2, 0\leq z \leq l \}.$$
Let $T$ be the solid torus by identifying $(x,y,0)$ and $(x,y, l)$
in $C$.  If $l$ is small enough depending only on $r$, then the
Hawking mass of $T$ is positive. Moreover,  if $T$ has an
admissible extension, then the Bartnik mass is also positive.
\end{exam}
To prove that $m_H(T)>0$, it is easy to see that the mean
curvature  $H$ of $\partial C$ is $\frac{1}{r}$, and
$$\int_{\partial C}H^2 =\frac{2\pi l}{r},$$
Suppose  $l<8r$, hence
$$\int_{\partial C}H^2 <16\pi. $$
We have $m_H(T)>0$.

Now suppose that $T$ has an admissible extension,  we may conclude
that $m_B(T)>0$ because of the following consequence of Theorem
\ref{bmasslowerbound}:

 \begin{prop}\label{torus} Suppose $(\Omega,g)$ has an admissible extension
  so that $\p\Omega$ is connected. Suppose $\p\Omega$ has positive mean
 curvature and suppose $m_H(\p\Omega)>0$. Then $m_B(\Omega)>0$.
\end{prop}
\begin{proof} Choose $\e_0>0$ be small enough such that
$\p\Omega_\e$ is smooth with positive mean curvature for all
$0<\e<\e_0$, where
$$
\Omega_\e=\{x\in \Omega|\ d(x,\p\Omega)>\e\}.
$$
Moreover, by choosing $\e_0$ small enough, we may also assume that
$m_H(\p\Omega_\e)>0$ for all $0<\e<\e_0$. Then each $\Omega_\e$ is
a minimizing hull in $\Omega$ because
 $\Omega\setminus\Omega_{\e_0}$ is foliated by positive mean
 curvature surfaces. Hence $m(\Omega_\e;\Omega)>0$ provided $\e>0$
 is small enough. In particular, $m(\Omega)>0$.
  By Theorem \ref{bmasslowerbound}, the proposition
 follows.
\end{proof}

Even though it is unclear if the solid torus $T$ in Example
\ref{cylinder} has an admissible extension. One can prove that $T$
has  no {\it smooth} static extension. Namely, there is no AF
manifold $(M,g)$ with one end satisfying the following:
\begin{itemize}
  \item [(a)] $g$ is smooth.
  \item [(b)] $T$ is isometrically embedded in $M$.
  \item [(c)] There is a smooth function $u$ defined on $M\setminus$interior
  of $T$ such that $u\to1$ at infinity, $Ric=u^{-1}\nabla^2u$ and $\Delta u=0$.
\end{itemize}
In fact, if $(M,g)$ is such an extension, then by \cite{M2},
$u=1-\frac m{r}+O(r^{-2})$ where $r$ is the Euclidean distance in
a coordinate system near infinity in the definition of AF manifold
and $m$ is the ADM mass of $M$. Moreover,  let $\Sigma=\p T$, then
 on $\Sigma$ (see   \cite{M2}):
\begin{equation}\label{torusboundary}
\Delta_\Sigma u+H \frac{\p u}{\p\nu}=0
\end{equation}
where $\Delta_\Sigma$ is the Laplacian on $\Sigma$, $\nu$ is the
unit outward normal and $H>0$ is the mean curvature of $\Sigma$.
Here we have used the fact that $\Sigma$ is flat, $H$ is constant
and $T$ is locally flat. Hence we have
\begin{equation*}
    \begin{split}
    0&=\int_{M\setminus T}\Delta u\\
    &=4\pi m-\int_{\Sigma}\frac{\p u}{\p \nu}\\
    &=4\pi m+H^{-1}\int_{\Sigma}\Delta_\Sigma u\\
    &=4\pi m.
    \end{split}
\end{equation*}
Hence $M$ must be $\R^3$. This is impossible by the uniqueness of
embedding.

\bibliographystyle{amsplain}

\end{document}